\documentclass{article}

\usepackage{arxiv}
\usepackage{newproof}

\usepackage[utf8]{inputenc} 
\usepackage[T1]{fontenc}    
\usepackage{hyperref}       
\usepackage{url}            
\usepackage{booktabs}       
\usepackage{amsfonts}       
\usepackage{nicefrac}       
\usepackage{microtype}      
\usepackage{lipsum}
\usepackage{amssymb}

\title{An elementary proof for dynamical scaling for certain fractional non-homogeneous Poisson processes}

\author{
  Markus Kreer \\ Feldbergschule, Oberh\"{o}chstadter Strasse 20, D-61440 Oberursel, Germany\\
  \texttt{mkreer@feldbergschule.de} \\
}

\newtheorem{thm}{Theorem}
\newtheorem{lem}[thm]{Lemma}
\newtheorem{cor}[thm]{Corollary}

\begin{document}

\maketitle

\begin{abstract}
Dynamical scaling is an asymptotic property typical for the dynamics of first-order phase transitions in physical systems and related to self-similarity. Based on the integral-representation for the marginal probabilities of a fractional non-homogeneous Poisson process introduced by Leonenko et al. (2017) and generalising the standard fractional Poisson process, we prove the dynamical scaling under fairly mild conditions. Our result also includes the special case of the standard fractional Poisson process.
\end{abstract}

\keywords{method of steepest decent  \and M-Wright function \and asymptotic analysis \and self-similarity}

\section{Preliminaries and main results}
Whereas the regular Poisson process  can be generalized in a straight-forward manner to a non-homogeneous Poisson process with a nicely behaved time-dependent intensity function by a mere time change using the chain rule, this trick does not work for fractional derivatives (e.g. \cite{Tarasov2016}). The original fractional Poisson process (fPp) as firstly introduced by \cite{Laskin2003} cannot be generalized using his approach. A possible way to go forth is following \cite{Scalas2017} who introduced  the marginals for a fractional non-homogeneus Poisson process (fnhPp) by equation Eq.(2.18) in their paper
\begin{eqnarray}
P_{\beta}(n,t) = \int_{0}^{\infty} du~ \frac{ e^{-\Lambda(u)} \Lambda(u)^n}{n!} h_{\beta}(t,u)
\label{eq:non-hom_fPp}
\end{eqnarray}
where $0<\beta<1$ and $\Lambda(t)=\int_{0}^{t} du~ \lambda(u)$ with an intensity function $\lambda(\cdot):[0,\infty)\rightarrow [0,\infty)$ and the density of the inverse stable subordinator (e.g. (\cite{Meerschaert2013})
\begin{eqnarray}
h_{\beta}(t,x)=\frac{t}{\beta x^{1+1/\beta}}~ g_{\beta}\left( \frac{t}{x^{1/\beta}}\right)
\label{eq:h_beta}
\end{eqnarray}
and
\begin{eqnarray}
g_{\beta}(u)=\frac{1}{\pi} \sum_{j=1}^{\infty} (-1)^{j+1}\frac{\Gamma(\beta j+1)}{\Gamma(j+1)}\frac{1}{u^{\beta j +1}}\sin(\beta \pi j)
\label{eq:g_beta}
\end{eqnarray}
The authors in \cite{Scalas2017} show in their Corollary 1 that the marginals Eq.(\ref{eq:non-hom_fPp}) satisfy the following system of fractional differential-integral equations
\begin{eqnarray}
\frac{d^{\beta}}{dt^{\beta}} P_{\beta}(n,t) & = &\int_{0}^{t} du ~\lambda(u) \left( -\frac{ e^{-\Lambda(u)} \Lambda(u)^n}{n!}+\frac{ e^{-\Lambda(u)} \Lambda(u)^{n-1}}{(n-1)!}\right) h_{\beta}(t,u),    
\hspace{1cm} n= 1,2,3,...
\nonumber
\\
\label{eq:nhfPp_eq}
\\
\frac{d^{\beta}}{dt^{\beta}} P_{\beta}(0,t)~ & = &~ -\int_{0}^{t} du ~\lambda(u) e^{-\Lambda(u)} h_{\beta}(t,u)    \hspace{5.4cm} n=0
\nonumber
\end{eqnarray}
with initial conditions $P_{\beta}(0,0)=1$ and $P_{\beta}(n,0)=0$ for $n=1,2,...$. The fractional derivative is the fractional Caputo-Djrbashian derivative which is defined as (e.g. \cite{Podlubni1999})
\begin{eqnarray}
\frac{d^{\beta}}{dt^{\beta}} F(t) = \frac{1}{\Gamma(1-\beta)} \int_0^t \frac{ds}{(t-s)^{\beta}}~\frac{d}{ds}f(s)
\nonumber
\end{eqnarray} 

As expected these equations Eq.(\ref{eq:nhfPp_eq}) reduce to the usual fractional Poisson process with Kolmogorov-Feller equations as discussed  by \cite{BeghinOrsingher2009} (see also \cite{Meerschaert2011} for a general theory, and \cite{Kreer2014} with an emphasis of their relation to Smoluchowsky-type coagulation-fragmentation equations) 
when $\lambda(t)$ is chosen to be just a positive constant $\lambda$,
\begin{eqnarray}
\frac{d^{\beta}}{dt^{\beta}} P_{\beta}(n,t) = \lambda \left( P_{\beta}(n-1,t)-P_{\beta}(n,t)\right)
\label{eq:fPp_eq}
\end{eqnarray}
where we have $0<\beta<1$ and $P_{\beta}(-1,t)=0$ for convenience and the same initial conditions as above. This fact proven by \cite{Scalas2017} will be demonstrated independently as we proceed with proving our main theorem.
\begin{thm} Let the function $\Lambda(\cdot)$ satisfy the following conditions 

(i) $\lim_{x\rightarrow\infty} \Lambda(x) =\infty$ and

(ii) There exists some positive constant $c$ such that as $x\rightarrow \infty $ 
\begin{eqnarray}
\frac{\Lambda'(x)}{\Lambda(x)} = \frac{c}{x} + \mathcal{O}\left(\frac{1}{x\log{x}}\right)   
\nonumber
\end{eqnarray}

Then for $0<\beta<1$ the marginal probabilities of the fnhPp as given in Eq.(\ref{eq:non-hom_fPp}) have the following limiting behaviour
\begin{eqnarray}
\lim_{n\rightarrow\infty, n=\Lambda(z_0 t^{\beta})} n P_{\beta}(n,t) =\frac{z_0}{c} M_{\beta}(z_0)
\nonumber
\end{eqnarray}
where for any positive $z_0>0$ we relate the positive $n$ to any positive $t$ by the equation $n=\Lambda(z_0 t^{\beta})$,  and $M_{\beta}(\cdot)$ is the M-Wright function or Mainardi function defined by
\begin{eqnarray}
M_{\beta}(z) =\frac{1}{\pi} \sum_{j=1}^{\infty} \frac{(-z)^{j-1}}{(j-1)!}\Gamma(\beta j)\sin{(\pi \beta j)}
\nonumber
\end{eqnarray}
\label{thm2}
\end{thm}

From Theorem \ref{thm2} we obtain immeadiately a useful 

\begin{cor} For any positive number $r$ let $\Lambda(x)= x^r$. Then  for $0<\beta<1$ the marginal probabilities as given in Eq.(\ref{eq:non-hom_fPp}) have the following limiting behaviour
\begin{eqnarray}
\lim_{t\rightarrow\infty, n t^{-r\beta}=z_0^r } t^{r\beta} P_{\beta}(n,t) =\frac{z_0^{1-r}}{r} M_{\beta}(z_0)
\nonumber
\end{eqnarray}
for any $z_0>0$ relating $t$ and $n$ by the condition $n t^{-r\beta}=z_0^{r}$. In particular for $r=1$ we obtain the dynamical scaling result for the standard  fPp satisfying the Kolmogorov-Feller equations Eq. (\ref{eq:fPp_eq}) with $\lambda=1$.
\label{cor3}
\end{cor}

Note that for $\beta=1$ and $\lambda(t)=1$, when we recover the standard Poisson process, the function $M_{\beta}(z)=0$ and our proof fails in this case.
\\

\section{Proof of Theorem \ref{thm2}}
We first show that the representation for the marginals as given in Eq.(\ref{eq:non-hom_fPp}) together with Eq.(\ref{eq:h_beta})-(\ref{eq:g_beta})   by \cite{Scalas2017}  can be expressed in terms of the M-Wright function (Mainardi function). Setting in Eq.(\ref{eq:non-hom_fPp}) for convenience $u=z t^{\beta}$ we get
\begin{eqnarray}
P_{\beta}(n,t) = t^{\beta}\int_{0}^{\infty} dz~ \frac{ e^{-\Lambda(z t^{\beta})} \Lambda(z t^{\beta})^n}{n!} h_{\beta}(t,z t^{\beta})
\nonumber
\end{eqnarray}
Because of  Eq.(\ref{eq:h_beta})-(\ref{eq:g_beta}) we obtain after some straight forward computations
\begin{eqnarray}
t^{\beta}  h_{\beta}(t,z t^{\beta}) 
& = &  \frac{1}{\beta z^{1+1/\beta}}~ g_{\beta}\left( \frac{1}{z^{1/\beta}}\right) \nonumber \\
& = & \frac{1}{\beta z^{1+1/\beta}}~ \frac{1}{\pi} \sum_{j=1}^{\infty} (-1)^{j+1}\frac{\Gamma(\beta j+1)}{\Gamma(j+1)} \frac{1}{\left(\frac{1}{z^{1/\beta}}\right)^{\beta j +1}}\sin(\beta \pi j) 
\nonumber \\
& = & \frac{1}{\pi} \sum_{j=1}^{\infty} (-1)^{j-1}\frac{\Gamma(\beta j)}{\Gamma(j)} z^{j-1} \sin(\beta \pi j)  =M_{\beta}(z)
\nonumber
\end{eqnarray}
The last equality follows from the series expansion of the M-Wright function (e.g. Eq. (10) in \cite{Mainardi2020})
Therefore we have as a backbone of our investigation the identity
\begin{eqnarray}
P_{\beta}(n,t) = \int_{0}^{\infty} dz~ \frac{ e^{-\Lambda(z t^{\beta})} \Lambda(z t^{\beta})^n}{n!} M_{\beta}(z)
\label{eq:non-hom_fPp1}
\end{eqnarray}
Note that when $\Lambda(\cdot)$ is a linear function we obtain a useful Lemma for fPp firstly stated by \cite{GorenfloMainardi2012}:
\begin{lem} (Gorenflo and Mainardi) For the fractional Poisson probabilities $P_{\beta}(n,t)$ we have
the integral representation (the subordination integral)
\begin{eqnarray}
P_{\beta}(n,t) & = &\frac{1}{n!}\int_{0}^{\infty} dt_{*} t_{*}^n ~ t^{-\beta} \exp{(-t_{*})} M_{\beta}(t_{*}t^{-\beta})
\nonumber \\
& = & \frac{1}{n!} t^{n \beta} \int_{0}^{\infty} dz ~ z^n e^{-z t^{\beta}} M_{\beta}(z)
\label{eq:int_repres}
\end{eqnarray}
where we have set for convenience $z=t_{*}t^{-\beta}$ as new integration variable.
\label{thm3}
\end{lem}

Thanks to Lemma 3 and its Eq.(\ref{eq:int_repres}), following as a special case from Eq.(\ref{eq:non-hom_fPp1}), we see immediately that without loss of generality if $\lambda(t)=1$ and consequently $\Lambda(t)=t$ is linear in $t$, the integration in the fractional differential-integral equations Eq.(\ref{eq:nhfPp_eq})  on the right-hand side can be carried out explicitly yielding the term $\left( P_{\beta}(n-1,t)-P_{\beta}(n,t)\right)$ and we have thus recovered the fractional Kolmogorov-Feller equations Eq.(\ref{eq:fPp_eq}) from \cite{BeghinOrsingher2009} as claimed. Therefore our dynamical scaling result covers also the standard fPp. We finally want to emphasize the fact that Lemma 3 is of highest practical relevance in computing the marginal probabilities of the fPp as functions of time for applications in statistics or physics.

{\em Proof of Theorem \ref{thm2}}

We write Eq.(\ref{eq:non-hom_fPp1}) as follows
\begin{eqnarray}
P_{\beta}(n,t) = \frac{1}{n!} \int_{0}^{\infty} dz~  e^{n f(z|t)} M_{\beta}(z)
\nonumber
\end{eqnarray}
Here we have defined an auxilary function by
\begin{eqnarray}
f(z|t)= \log{\Lambda(zt^{\beta})}-\frac{1}{n} \Lambda(z t^{\beta})
\label{eq:f_funct}
\end{eqnarray}
and we treat $t$ as a parameter.
We compute now the first and second derivative of $f(\cdot|t)$ with respect to the first argument, namely
\begin{eqnarray}
f'(z|t) 
= \Lambda'(z t^{\beta}) t^{\beta}\left[\frac{1}{\Lambda(z t^{\beta})}-\frac{1}{n}\right]
\label{eq:f_deriv1}
\end{eqnarray}
and
\begin{eqnarray}
f''(z|t) = \Lambda''(z t^{\beta}) t^{2\beta}\left[\frac{1}{\Lambda(z t^{\beta})}-\frac{1}{n}\right]
-\left( \frac{\Lambda'(z t^{\beta})t^{\beta}}{\Lambda(z t^{\beta})} \right)^2
\label{eq:f_deriv2}
\end{eqnarray}
Due to the properties of $\Lambda(\cdot)$, condition (i), there exists for any positive $n$ and for any $z_0>0$ a positive number $t$ solving , $f'(z_0|t)=0$, yielding from Eq.(\ref{eq:f_deriv1})
\begin{eqnarray}
n = \Lambda(z_0 t^{\beta})
\label{eq:n_Lambda}
\end{eqnarray}
Note that from Eq.(\ref{eq:n_Lambda}) we have with Eq. (\ref{eq:f_funct}) that  $f(z_0|t)=\log{\Lambda(z_0 t^{\beta})}-\Lambda(z_0 t^{\Lambda})/n = \log{n} -1$.
And  we see from Eq.(\ref{eq:f_deriv2}) that $f''(z_0|t)<0$ and thus $z_0>0$ is a positive maximum of $f(\cdot|t)$, i.e. using condition (ii) on $\Lambda(\cdot)$ we get
\begin{eqnarray}
\sqrt{-f''(z_0|t)} = \frac{\Lambda'(z_0 t^{\beta}) t^{\beta}}{\Lambda(z_0 t^{\beta})} = \frac{c}{z_0 t^{\beta}} t^{\beta} + t^{\beta}\mathcal{O}\left(\frac{1}{t^{\beta}\log{t}}\right) 
\end{eqnarray}

Because the M-Wright function $M_{\beta}(\cdot)$ is analytic on the real line (e.g. \cite{GorenfloLuchenkoMainardi1999}, \cite{Podlubni1999}) we can evaluate the integral using the method of steepest decent because the conditions for applying \cite{Small2010} Chapter 6 Proposition 2 and Proposition 3 are satisfied:
\begin{eqnarray}
\int_0^{\infty} dz ~ e^{n f(z|t) }M_{\beta}(z) & = & \frac{\exp{(n f(z_0|t))}}{\sqrt{n}} \left(M_{\beta}(z_0)~\sqrt{\frac{2\pi}{-f''(z_0|t)}} +\mathcal{O}\left(\frac{1}{n}\right)\right)
\nonumber \\
& = & 
\frac{e^{n f(z_0|t)}}{\sqrt{n}}\sqrt{2\pi} \left(M_{\beta}(z_0)~\frac{z_0}{c} +\mathcal{O}\left(\frac{1}{n}\right)+\mathcal{O}\left(\frac{1}{\log{t}}\right)\right)
\nonumber
\end{eqnarray}
and with Stirling's formula 
\begin{eqnarray}
n! = \sqrt{2 \pi n}  \left(\frac{n}{e}\right)^n \left(1+\mathcal{O}\left(\frac{1}{n}\right)\right)
\nonumber
\end{eqnarray}
 and Eq.(\ref{eq:n_Lambda}) we obtain altogether
\begin{eqnarray}
 \frac{1}{n!} \int_0^{\infty} dz ~ e^{n f(z|t) }M_{\beta}(z) 
& = & \frac{1}{\sqrt{2 \pi n}}  \left(\frac{n}{e}\right)^{-n}\left(1+\mathcal{O}\left(\frac{1}{n}\right)\right)  \cdot
\nonumber \\
&  & \cdot
e^{n f(z_0|t)} \sqrt{\frac{2\pi}{n}} \left(M_{\beta}(z_0)~\frac{z_0}{c} +\mathcal{O}\left(\frac{1}{n}\right)+\mathcal{O}\left(\frac{1}{\log{t}}\right)\right)
\nonumber \\
& = & \frac{1}{n}
\left(M_{\beta}(z_0)~\frac{z_0}{c} +\mathcal{O}\left(\frac{1}{n}\right)+\mathcal{O}\left(\frac{1}{\log{t}}\right)\right)
\nonumber
\end{eqnarray}
where we have used  $f(z_0|t)= \log{n} -1$ in the last line. Multiply both sides with $n$ and note that for large increasing $n$ by the properties of $\Lambda(\cdot)$ also $t$ increases and $z_0>0$ remains a maximum as long as $n$ and $t$ are related by Eq. (\ref{eq:n_Lambda}). This concludes the proof. 
$\blacksquare$


\section{Outlook and application}
Since \cite{ViscekFamily1984} the dynamical scaling hypothesis became a ``regular'' feature in the physical theory of  phase transitions. It was originally discovered in numerical studies and not long afterwards a mathematical proof for the simple Smoluchowski coagulation equations with constant kernel was given by \cite{KreerPenrose1994}. The dynamical scaling hypothesis is expected to hold for generalized Smoluchowski equations but not proven yet. The work of \cite{Kreer2014} related the fPp to a certain class of Smoluchowski equations and thus a dynamical scaling of the marginal probabilities would come to no surprise. Our current work has answered this in the affirmative even for the more general class of fnhPp as introduced by \cite{Scalas2017}. This generalisation has the advantage that the dynamical scaling exponent of time $t$ is not limited to positive numbers between 0 and 1 (which is the range the order of the fractional derivative can take) but can be any positive number due to the fact that the the choice of the intensity function $\Lambda(\cdot)$ allows for some additional freedom.

How would one apply a fnhPp to the problem of phase transitons? In our case  we may  imagine a solid crystal consisting of single monomers bound in a regular lattice structure. When the crystal is heated the lattice takes up heat and the monomeres start to oscillate and finally break free. At the end, the lattice has been dissolved because all monomers have gone. The very simple description of this ``phase transition'', in which the ordered state (= crystall lattice) has undergone a transition to an unordered state (free-moving monomers corresponding to a liquid or a gas) is to count the number of monomers breaking free from the remaining crystal structure, one at a time. The reverse of this model would be monomers condensing to a newly formed crystal. In fact, scientists have used this simple counting model to describe the nucleation of lysozyme and paracetamol in aqueous droplets:  \cite{Goh2010} describe the nucleation by a Poisson process and estimate its parameters. 

A simple computation shows that a standard Poisson process does not possess a regular dynamical scaling limit, namely with $\beta=1$ in Eq.(\ref{eq:int_repres}) we have after some straight-forward computations
\begin{eqnarray}
\lim_{t\rightarrow\infty, t^{-1}n=z_{0}} t^{1/2} P_{1}(n,t) = 
\left\{
\begin{array}{ll}
\frac{1}{\sqrt{2\pi}} \hspace{0.5cm} \textrm{for} \hspace{0.2cm}  z_{0}=1\\
0 \hspace{1cm} \textrm{otherwise} 
\end{array}
\right.
\nonumber
\end{eqnarray}
With respect to the dynamical scaling hypothesis for phase transitions it would be very interesting to analyse the nucleation data of \cite{Goh2010} using a fractional Poisson process as model.

\section*{Acknowledgements}

Stimulating discussions with Lukas Kreer and Theo Manoussos (both Johannes-Gutenberg Universit\"{a}t, Mainz) and Dr. Ay\c{s}e K{\i}z{\i}lers\"{u} (University of Adelaide) are gratefully acknowledged.

This research did not receive any specific grant from funding agencies in the public, commercial or not-for-profit sectors.

\bibliographystyle{unsrt}  

\end{document}